\documentclass{article}

\usepackage{delarray,verbatim,enumerate,a4wide}
\usepackage{amsmath,amsthm,amstext,amsbsy,amssymb,amsfonts,amscd}
\usepackage[all]{xy}
\usepackage[frenchb]{babel}
\usepackage[T1]{fontenc}
\usepackage[utf8]{inputenc}

\usepackage{scalerel}

\usepackage{upgreek,bbm}

\usepackage[small]{titlesec}

\usepackage{color}
\definecolor{vert}{rgb}{0.1,0.4,0.2}
\usepackage[colorlinks=true,linkcolor=blue,citecolor=vert]{hyperref}

\usepackage{calligra}
\DeclareFontShape{T1}{calligra}{m}{n}{<->s*[0.95]callig15}{}
\DeclareMathAlphabet{\mathscr}{T1}{calligra}{m}{n}

\newtheorem{Th}{Théorème}[]

\newtheorem{Prop}[Th]{Proposition}
\newtheorem{Cor}[Th]{Corollaire}

\newtheorem{Def} [Th]{Définition}

\newtheorem*{ThP}{Théorème Principal}
\newtheorem*{CoP}{Corollaire}
\newtheorem*{Def*}{Définition}

\def\Preuve{\noindent {\it Preuve.~}}

\def\Remarque{\smallskip\noindent {\it Remarque.~}}

\def\Nota{\smallskip\noindent {\it Nota.~}}

		\def\QQ{\mathbb Q}	
\def\NN{\mathbb N}	\def\ZZ{\mathbb Z}		
\def\F2{\mathbb{F}_2}	\def\Z2{\mathbb{Z}_2}	\def\lT{\ov{\mathbb T}_\ell}	
\def\Zl{{\mathbb{Z}_\ell}} 	\def\Ql{\mathbb{Q}_\ell}	\def\Tl{\mathbb{T}_\ell}	\def\Fl{\mathbb{F}_\ell}	\def\hSl{\hat{\mathcal{S}\ell}}

 				\def\U{\mathcal  U}	\def\F{\mathcal  F}	
\def\J{\mathcal  J}  		\def\C{\mathcal  C}		\def\R{\mathcal  R}		\def\H{\mathcal  H}
 	\def\Pl{\mathcal  P\ell}  	\def\Cl{\mathcal  C\!\ell}	
		\def\T{\mathcal  T}			

		\def\mmu{\boldsymbol{\upmu}}

		\def\p{{\mathfrak p}}		\def\x{{\mathfrak x}}

\def\W{W\!K}			

\def\wi{\widetilde}	\def\ov{\overline}	
	\def\div{\operatorname{div}}
	\def\deg{\operatorname{deg}}		
\def\Gal{\operatorname{Gal}}	\def\Log{\operatorname{Log}}		\def\Rad{\operatorname{Rad}}
		\def\Hom{\operatorname{Hom}}

\newcommand\scale[2]{\vstretch{#1}{\hstretch{#1}{#2}}}

\newcommand\si[1]{\scale{.7}{#1}}
\newcommand\ph{{\phantom{*}}}
\newcommand\ab{{\scale{.8}{\rm ab}}}	\newcommand\lc{{\scale{.8}{\rm lc}}}	
		\newcommand\bp{{\scale{.8}{\rm bp}}}
\newcommand\s{{\scale{.7}{\mathcal S}}}		\newcommand\lz{{\scale{.8}{\rm lz}}}
\def\%{{\scale{.8}{\infty}}}			

\makeatletter
\newcommand*\wt[2][0.2ex]{%
        \begingroup
        \mathchoice{\wt@helper{#1}{#2}{\displaystyle}{\textfont}}
                   {\wt@helper{#1}{#2}{\textstyle}{\textfont}}
                   {\wt@helper{#1}{#2}{\scriptstyle}{\scriptfont}}
                   {\wt@helper{#1}{#2}{\scriptscriptstyle}{\scriptscriptfont}}%
        \endgroup
        #2%
}
\newcommand*\wt@helper[4]{%
        \def\currentfont{\the#41}%
        \def\currentskewchar{\char\the\skewchar\currentfont}%
        \setbox\tw@\hbox{\currentfont$#2$\currentskewchar}%
        \dimen@ii\wd\tw@
        \setbox\tw@\hbox{\currentfont$#2${}\currentskewchar}%
        \advance\dimen@ii-\wd\tw@
        \rlap{\raisebox{-#1}{$\m@th#3\kern\dimen@ii\widetilde{\phantom{#2}}$}}%
}
\makeatother

\def\wE{\,\wt[0.1ex]{\!\mathcal E}}		\def\wc{\wt[0.3ex]{\mathfrak C}}
\def\wU{\wt[0.2ex]{\mathcal U}}	\def\wC{\wt[0.1ex]{\mathcal C}}
\def\wJ{\,\wt[0.2ex]{\!\mathcal J}}	\def\wCl{\wt[0.1ex]{\mathcal C\!\ell}} \def\wDl{\wt[0.2ex]{\mathcal D\!\ell}}

\begin{document}

\title{\Large\bf Annulateurs de Stickelberger des groupes de classes logarithmiques}

\author{ Jean-François {\sc Jaulent} }
\date{}
\maketitle
\bigskip

{\small
\noindent{\bf Résumé.} Étant donnés un nombre premier impair $\ell$ et un corps de nombres abélien $F$ contenant les racines $\ell$-ièmes de l'unité, nous montrons que l'idéal de Stickelberger annule la composante imaginaire du $\ell$-groupe des classes logarithmiques et que son reflet annule la composante réelle du module de Bertrandias-Payan. Nous en tirons une preuve très simple d'un théorème d'annulation pour les $\ell$-noyaux étales sauvages.

\

\noindent{\bf Abstract.} For any odd prime number $\ell$ and any abelian number field $F$ containing the $\ell$-th roots of unity, we show that the Stickelberger ideal annihilates the imaginary component of the $\ell$-group of logarithmic classes and that its reflection annihilates the real component of the Bertrandias-Payan module. This leads to a very simple proof of annihilation results for the so-called wild {\it étale} $\ell$-kernels of $F$.

\


}



\section*{Introduction}
\addcontentsline{toc}{section}{Introduction}

Sous sa forme originelle, le classique théorème de Stickelberger (cf. e.g. \cite{Wa}) affirme que pour tout corps abélien $F_{\si{0}}$ de groupe de Galois $G_{F_{\si{0}}}$, l'idéal $\mathcal S_{F_{\si{0}}}$ de l'algèbre $\ZZ[G_{F_{\si{0}}}]$ engendré par les éléments de Stickelberger tordus annule le groupe des classes d'idéaux de $F$. Lorsque $F$ est imaginaire, il est bien connu que ces éléments le sont aussi, en ce sens qu'ils se factorisent par $1-\bar\tau$, où $\bar\tau$ désigne la conjugaison complexe. Il suit de là que ce théorème n'apporte aucune information {\em directe} sur le groupe des classes du sous-corps réel $F^+$, lequel est banalement tué par $1-\bar\tau$.

Pour aller plus loin, il est utile de fixer un nombre premier $\ell$ (que nous prendrons ici impair pour éviter les complications techniques liée au cas $\ell=2$) et de se concentrer sur les $\ell$-parties des groupes considérés. Dans cette perspective, il est naturel de raisonner en présence du groupe $\mmu_\ell$ des racines $\ell$-ièmes de l'unité en remplaçant le corps de départ par le sur-corps $F=F_{\si{0}}[\zeta_\ell]$; auquel cas les $\Zl$-invariants de $F_{\si{0}}$ se déduisent de ceux de $F$ par redescente via l'idempotent $e_{F/F_{\si{0}}} = \frac{1}{[F:F_{\si{0}}]}\sum_{\sigma\in\Gal(F/F_{\si{0}})}\sigma$. Et ce n'est pas restreindre la généralité que travailler sur $F$.\smallskip

Partons donc d'un corps abélien $F$ contenant les racines $\ell$-ièmes de l'unité; notons $e_\pm=\frac{1}{2}(1\pm\bar\tau)$ les idempotents associés à la conjugaison complexe $\bar\tau\in G_F=\Gal(F/\QQ)$; et convenons de dire qu'un $\Zl[G_F]$-module est réel lorsqu'il est fixé par $e_+$, imaginaire lorsqu'il l'est par $e_-$. Cela étant:\smallskip

-- D'un côté, le théorème de Stickelberger affirme que l'idéal imaginaire $\mathcal S_F$ de l'algèbre $\ZZ[G_F]$ annule la composante imaginaire du $\ell$-groupe des classes d'idéaux $\,\Cl_F$ de $F$; autrement dit, que pour tout $m$ assez grand, sa réduction modulo $\ell^m$ annule la composante imaginaire de $\,\Cl_F$.\smallskip

-- D'un autre côté, Gras et Oriat (cf. \cite{Gra1,Or1,Gra5}) ont montré que le reflet modulo $\ell^m$ de l'idéal $\mathcal S_F$ est un idéal réel de l'algèbre quotient  $\ZZ/\ell^m\ZZ[G_F]$ qui annule pour tout $m$ assez grand la composante réelle du sous-groupe de torsion $\,\T^\%_F$ du pro-$\ell$-groupe des classes infinitésimales de $F$.\smallskip

Dans ces deux formulations le $\ell$-groupe de classes $\,\Cl_F$ d'une part et $\ell$-groupe $\,\T^\%_F$ d'autre part apparaissent comme en miroir, alors même qu'il n'existe pas de dualité naturelle entre eux : en particulier, ils n'ont généralement ni même ordre, ni même rang.\medskip

Le but de la présente note est de présenter une formulation duale originale de ces théorèmes d'annulations, obtenue en remplaçant $\,\Cl_F$ par son analogue logarithmique $\,\wCl_F$ et $\,\T^\%_F$ par ce qu'il est convenu d'appeler le module de Bertrandias-Payan $\,\T^\bp_F$ (cf. e.g. \cite{Gra3,J23,Ng2}).\smallskip

Nous énonçons d'abord, pour un corps de nombres arbitraire $F$, les résultats bien connus sur les groupes $\,\wCl_F$ dont nous avons besoin par la suite, ainsi que quelques autres sur les groupes $\T^\bp_F$ dans le contexte kummérien, issus essentiellement de \cite{Gra2}, \cite{J23} et \cite{J28}. Puis, pour $\ell$ impair et $F$ abélien contenant $\mmu_\ell$, donc de conducteur $f_F$ multiple de $\ell$, nous reprenons la définition standard des éléments de Stickelberger tordus $\s_F^{\,c}$, telle qu'énoncée par exemple dans \cite{Gra5}, et nous définissons leurs reflets $\s_F^{\,c*}$ dans l'involution du miroir par montée et descente dans la $\Zl$-tour cyclotomique. Procédant ensuite par induction à partir du résultat classique pour la partie imaginaire, puis par dualité de Kummer pour la partie réelle, nous obtenons alors le Théorème Principal suivant:

\begin{ThP}\label{TP}
Soient $\ell$ un nombre premier impair et $F$ un corps abélien contenant $\mmu_\ell$. Notons $\,G_F$ son groupe de Galois, $F_\%=\bigcup F_n$ sa $\Zl$-extension cyclotomique, $F_n^\bp$ le compositum des $\ell$-extensions cycliques de $F_n$ localement $\Zl$-plongeables, $ F_n^z$ le compositum des $\Zl$-extensions et  $F_n^\lc$ la plus grande pro-$\ell$-extension abélienne de $F_n$ complètement décomposée partout sur $F_\%$.\smallskip

Notons $F_\%^\lc=\bigcup F_n^\lc$, puis $\,\wC_{F_\%}=\Gal(F_\%^\lc/F_\%)$ le module de Kuz'min-Tate, $\,\wc_{F_\%}=\Rad(F_\%^\lc/F_\%)$ son dual de Kummer $\Hom(\,\wC_{F_\%},\mmu_{\ell^\%})$; écrivons $\,\Gamma_{\!n}=\Gal(F_\%/F_n)$ et $G_{F_n}=\Gal(F_n/\QQ)$.\smallskip


Définissons le symétrisé $\ell$-adique $\hat{\mathcal S}\ell_{F_n}$ de l'idéal de Stickelberger comme l'idéal de l'algèbre $\Zl[G_{F_n}]$ engendré par les éléments de Stickelberger tordus $\s_{F_n}^{\,c}$  et leurs reflets $\s_{F_n}^{\,c*}$ pour $c\in\ZZ_\ell^\times$.\smallskip
\begin{itemize}
\item[(i)] Le $\ell$-groupe des classes logarithmiques $\,\wCl_{F_n} \simeq \Gal(F_n^\lc/F_\%)$ s'identifie au quotient $^{\Gamma_{\!n}}\wC_{F_\%}$.\smallskip

\item[(ii)] Le module de Bertrandias-Payan $\,\T^\bp_{F_n} \simeq \Gal(F_n^\bp/F_n^z)$ au sous-groupe invariant $\wc_{F_\%}^{\,\Gamma_{\!n}}$.\smallskip

\item[(iii)] Et l'idéal  symétrisé $\hat{\mathcal S}\ell_{F_n}$, d'indice fini dans $\Zl[G_{F_n}]$, annule à la fois  $\,\wCl_{F_n}$ et $\,\T^\bp_{F_n}$.

\end{itemize}
\end{ThP}

\Nota Il résulte des calculs de Sinnott \cite{Si} que l'idéal  $\hat{\mathcal S}\ell_{F_n}$ est d'indice fini dans l'algèbre $\Zl[G_{F_n}]$.\medskip

Du fait des relations dans la $\Zl$-tour cyclotomique entre $\ell$-groupes de classes logarithmiques et $\ell$-noyaux étales supérieurs, qui reposent sur la description, déjà exploitée dans \cite{J44}, obtenue par Schneider dans \cite{Sc} en termes de théorie d'Iwasawa, nous en déduisons alors une preuve très simple d'un théorème d'annulation pour ces mêmes noyaux. Le Corollaire \ref{Noyaux étales} infra donne ainsi:

\begin{CoP}
Pour chaque $i\in\ZZ$, définissons le  $(-i)$-ième tordu à la Tate $\hSl_{F_n}^{(i)}$ du symétrisé $\ell$-adique $\hSl_{F_n}$ de l'idéal de Stickelberger  comme étant  l'idéal de l'algèbre $\Zl[G_{F_n}]$ engendré par les tordus à la Tate des éléments de Stickelberger $\s_{F_n}^{\,c}$ et de leurs reflet $\s_{F_n}^{\,c*}$ pour $c\in\ZZ_\ell^\times$.\par
L'idéal obtenu $\hSl_{F_n}^{(i)}$ annule le $i$-ième $\ell$-noyau étale sauvage $\W_{2i}({F_n})$ attaché au corps $F_n$.
\end{CoP}

Rappelons que, suite aux travaux de Soulé, Rost, Voevodsky et Weibel (cf. e.g. \cite{Sou, We}), les $\ell$-sous-groupes de Sylow des groupes de $K$-théorie $K_{2i}(\mathcal O_F)$ attachés à l'anneau des entiers d'un corps de nombres pour $i>0$, s'identifient à des groupes de cohomologie étale, ce qui permet de regarder les $\ell$-noyaux sauvages de Banaszak \cite{Ba} comme des noyaux étales supérieurs, suivant la terminologie introduite par Nguyen Quang Do \cite{Ng1} et d'étendre leur définition à tout $i\in\ZZ$.\smallskip

Le Corollaire redonne donc le résultat d'annulation établi par Coates et Sinnott \cite{CS} pour la composante imaginaire $\mathcal S\ell_F $ de $\hat{\mathcal S}\ell_F$ et $i>0$, mais vaut encore pour le symétrisé $\hat{\mathcal S}\ell_F$ et même pour tout $i\in\ZZ$, indépendamment des conjectures de finitude de Schneider.\smallskip

Précisons ici que c'est volontairement que nous nous sommes restreints dans cette note aux cas des corps absolument abéliens, pour lesquels la définition des annulateurs de Stickelberger est complètement élémentaire. On peut néanmoins imaginer que l'approche logarithmique que nous avons utilisée pourrait aussi bien être développée dans un cadre plus général, dans le contexte de ce qu'il est convenu d'appeler la conjecture équivariante, dans lequel plusieurs résultats profonds viennent d'être récemment obtenus (cf. \cite{DK,JN,Ni}).

Signalons enfin que, pour les corps abéliens réels, les sous-groupes d'unités circulaires logarithmiques, ainsi que ceux de normes cyclotomiques universelles, fournissent également des annulateurs pour les groupes de classes logarithmiques (cf. \cite{J66}), qui sont le pendant logarithmique des résultats de All \cite{All} prouvant et  généralisant la conjecture de Solomon \cite{Sol}. Il sera naturellement instructif d'illustrer numériquement les résultats obtenus pour les classes logarithmiques par ces deux approches, à l'instar de ce qui est fait par Gras \cite{Gra5} dans le cadre classique.


\newpage
\section{Présentation des classes et unités logarithmiques}

Soient $\ell$ un nombre premier donné et $F$ un corps de nombres.  À chaque place finie  $\p$ de $F$, il est attaché dans \cite{J28} une application à valeurs dans $\Zl$, définie sur le groupe multiplicatif $F_\p^\times$ du complété de $F$ en $\p$ par la formule:

\centerline{$\tilde\nu_\p (x_\p)\, =\,\nu_\p (x_\p)$, pour $\p\nmid\ell$;\quad et \quad $\tilde\nu_\p (x_\p)\, =  -\frac{1}{\deg\, \p}\,\Log_\ell (N_{F_\p/\Ql}(x_\p))$, pour $\p|\ell$;}\smallskip

\noindent où $\Log_\ell$ désigne le logarithme d'Iwasawa et $\deg\p$ est un facteur de normalisation calibré pour assurer que l'image de $F_\p^\times$ soit dense dans $\Zl$. Cette application induit un morphisme surjectif du {\em compactifié $\ell$-adique} du groupe multiplicatif $F_\p^\times$\smallskip

\centerline{$\R_{F_\p}=\varprojlim  F_\p^\times/F_\p^{\times\ell^n}$}\smallskip

\noindent sur $\Zl$ dont le noyau, dit {\em sous-groupe des unités logarithmiques} de $\R_{F_\p}$,\smallskip

\centerline{$\wU_{F_\p}\,=\,\{x_\p\in\R_{F_\p}\,|\, \tilde\nu_\p(x_\p)=0\}$}\smallskip

\noindent  s'identifie par la Théorie $\ell$-adique locale du corps de classes (cf. \cite{J31}) au sous groupe normique de $\R_{F_\p}$ associé à la $\Zl$-extension cyclotomique $F_\p^{\,c}$ de $F_\p$.
C'est donc l'analogue du groupe\smallskip

\centerline{$\U_{F_\p}\,=\,\{x_\p\in\R_{F_\p}\,|\, \nu_\p(x_\p)=0\}$}\smallskip

\noindent des unités de $\R_{F_\p}$, qui correspond, lui, à la $\Zl$-extension non-ramifiée de $F_\p$.
\medskip

Soit maintenant $\J _F$ le {\em $\ell$-adifié  du groupe  des idèles de $F$}, i.e. le produit $\J_F=\prod_\p ^{\si{\rm res}}\R_{F_\p}$ des compactifés  $\R _{F_\p}$ des groupes multiplicatifs des complétés $F_\p$, restreint aux familles $(x_\p)_\p$ dont presque tous les éléments tombent dans le sous-groupe unité $\,\U_F=\prod_\p \, \U _{F_\p}$.
La Théorie $\ell$-adique globale du corps de classes établit un isomorphisme de groupes topologiques entre le $\ell$-groupe des classes d'idèles de $F$, défini comme le quotient {$\J_F/\R_F$} de $\J_F$ par son sous-groupe principal $\R_F=\Zl\otimes_\ZZ F^\times$, et le groupe de Galois $\Gal(F^{\ab}/F)$ de sa pro-$\ell$-extension abélienne maximale $F^\ab$.\smallskip

Dans la correspondance ainsi établie (cf. \cite{J28,J31}),\smallskip

\begin{itemize}
\item[(i)] le groupe de normes associé à la $\Zl$-extension cyclotomique $F^c=F_{\%}$ de $F$ est le sous-groupe des idèles de degré nul: $\wJ_F=\{\x=(x_\p)_\p\in\J_F\,|\, \deg(\x)=\sum_\p\wi\nu_\p(x_\p)\deg\,\p=0\}$;\smallskip

\item[(ii)] le groupe de normes associé à la plus grande sous-extension $F^{\lc}$ de $F^{\ab}$ qui est localement cyclotomique (i.e. complètement décomposée sur $F_\%$ en chacune de ses places) est le produit $\,\wU_F\R_F$ du sous-groupe $\,\wU_F=\prod_\p \, \wU _{F_\p}$ des unités logarithmiques locales et de $\R_F$;\smallskip

\item[(iii)] le groupe de Galois $\Gal(F^{lc}/F_\%)$ s'identifie ainsi au quotient $\,\wCl_F= \wJ_F/\wU_F\R_F$ du groupe $\wDl _F = \wJ _F /\wU_F$ des diviseurs logarithmiques de degré nul par son sous-groupe principal $\Pl_F=\R_F\wU_F/\wU_F$, image canonique de $\R_F$;\smallskip

\item[(iv)] et le noyau $\wE _F=\R_F \cap\, \wU_F$ du morphisme $\wi\div:\;x\mapsto \sum_\p\wi\nu_\p(x)\,\p$ de $\R_F$ dans $\wDl_F$ est le sous-groupe des normes cyclotomiques (locales comme globales): $\wE_F=\bigcap_n N_{F_n/F}(\R_{F_n})$, où $F_n$ décrit les étages finis de la $\Zl$-tour cyclotomique $F_\%/F$.
\end{itemize}
\smallskip

\begin{Def*}
On dit que $\,\wCl_F\simeq \Gal(F^\lc/F_\%)$ est le $\ell$-groupe des classes logarithmiques du corps $F$ et que $\,\wE_F$ est le pro-$\ell$-groupe des unités logarithmiques globales.
\end{Def*}

Comme expliqué dans \cite{J28,J55}, la {\em conjecture de Gross-Kuz'min} (pour le corps $F$ et le premier $\ell$) revient à postuler la finitude de $\,\wCl_F$ ou, de façon équivalente, que le $\Zl$-rang de $\,\wE_k$ est le somme $r_{\si{F}}+c_{\si{F}}$ des nombres de places réelles et complexes de $F$; ce qui est le cas dès lors que $F$ est abélien.\smallskip

Enfin, du point de vue de la théorie d'Iwasawa, le groupe $\,\wCl_F$ s'interprète comme le quotient des genres ${}^\Gamma\C'_F$, relativement au groupe procyclique $\Gamma=\Gal(F_\%/F)$, du {\em module de Kuz'min-Tate}\smallskip

\centerline{$\C'_{F_\%} = \varprojlim \,\Cl'_{F_n}$}\smallskip

\noindent (noté $\,\T_F$ dans \cite{J55}), limite projective des $\ell$-groupes de $\ell$-classes d'idéaux attachés aux étages finis $F_n$ de la tour $F_\%/K$. 
En particulier on a aussi bien: $\,\C'_{F_\%}=\wC_{F_\%}$, où\smallskip

\centerline{$\wC_{F_\%} = \varprojlim \,\wCl_{F_n}$}\smallskip

\noindent est la limite projective des groupes logarithmiques  $\,\wCl_{F_n}={}^{\Gamma_{\!n}}\C'_{F_\%}$ avec $\Gamma_{\!n}=\Gamma^{\ell^n}=\Gal(F_\%/F_n)$.

\newpage
\section{Lien avec le module de Bertrandias-Payan}

Il est traditionnel depuis \cite{BP} d'appeler pro-$\ell$-extension de Bertrandias-Payan attaché à un corps de nombres $F$  le compositum $F^\bp$ des $\ell$-extensions cycliques $L$ de $F$ qui sont {\em localement} $\Zl$-plongeables. En d'autres termes, si $F^z$ désigne, lui, le compositum des $\Zl$-extensions de $F$, le corps $F^\bp$ est la plus grande pro-$\ell$-extension abélienne $F^\lz$ de $F$ qui est complètement décomposée sur $F^z$ en chacune de ses places (cf. \cite{J23}, 2\S b). 

\begin{Def*}
Le sous-groupe de torsion  $\T_F^\bp=\Gal(F^\bp/F^z)$ du groupe abélien $\H_F=\Gal(F^\bp/F)$ est un groupe fini, dit module de Bertrandias Payan.
\end{Def*}

Lorsque le corps de base $F$ contient les racines $2\ell$-ièmes de l'unité, le module $\,\T_F$ est quasiment en dualité kummérienne avec le groupe des classes logarithmiques $\,\Cl_F$. Expliquons cela:

\begin{itemize}
\item Pour chaque étage fini $F_n$ de la tour $F_\%/F$, le groupe de Galois $\H_{F_n}=\Gal(F_n^\bp/F_n)$ attaché au $n$-ième étage $F_n$ de la tour $F_\%/F$ est donné par l'isomorphisme:\smallskip

\centerline{$\H_{F_n}\simeq\J_{F_n}/\big(\prod_{\p^\ph_n}\mu_{\p_n}\big)\R_{F_n}$,}\smallskip

\noindent où $\p_n$ parcourt les places finies de $F_n$ et $\mu_{\p_n}$ désigne le $\ell$- groupe des racines de l'unité du complété $F_{\p_n}$ de $F_n$ en $\p_n$ (cf. \cite{J23}, 2\S b).

\item Introduisant le pseudo-radical $\,\mathfrak R_{F_n}=(\Ql/\Zl)\otimes_\ZZ F_n^\times= (\Ql/\Zl)\otimes_\Zl \R_{F_n}$, on voit alors que le module de Bertrandias-Payan $\,\T_{F_n}^\bp$ s'identifie au sous-groupe $\,\wc_{F_n}$ de $\,\mathfrak R_{F_n}$ défini par:\smallskip

\centerline{$\wc_{F_n} = \{ \ell^{-k}\otimes x_n \in \mathfrak R_{F_n} \;|\; x_n\in\big(\prod_{\p^\ph_n}\mu_{\p_n}\big)\J_{F_n}^{\ell^k}\}$}\smallskip

\item Et, par la théorie de Kummer, le corps $F_n$ contenant par hypothèse les racines $\ell^n$-ièmes de l'unité, le sous-groupe de $\ell^n$-torsion ${}_{\ell^n}\wc_{F_n}$ de $\,\wc_{F_n}$, s'interprète comme radical kummérien attaché à la sous-extension d'exposant $\ell^n$ de $F_n^\lc$ (cf. \cite{J23}, Déf. 2.6).
\end{itemize}\smallskip

Plus précisément ces pseudo-radicaux vérifient la théorie de Galois dans la tour $F_\%/F$, en ce sens que, posant $\,\mathfrak R_{F_\%}=(\Ql/\Zl)\otimes_\ZZ F_\%^\times= (\Ql/\Zl)\otimes_\Zl \R_{F_\%}$ puis\smallskip

\centerline{$\wc_{F_\%} = \{ \ell^{-k}\otimes x_\% \in \mathfrak R_{F_\%} \;|\; x_\%\in\big(\prod_{\p^\ph_\%}\mu_{\p_\%}\big)\J_{F_\%}^{\ell^k}\}$,}\smallskip

\noindent on a canoniquement: $\,\mathfrak R_{F_n}=\mathfrak R_{F_\%}^{\,\Gamma_{\!n}}$ et $\,\wc_{F_n}=\wc_{F_\%}^{\,\Gamma_{\!n}}$, pour $\Gamma_{\!n}=\Gal(F_\%/F_n)$ (cf. \cite{J23}, Prop. 1.2).

Maintenant, $\,\wc_{F_\%}$ est précisément le radical $\Rad(F_\%^\lc/F_\%)$ attaché à la pro-$\ell$-extension abélienne maximale de $F_\%$ qui est complètement décomposée en chacune de ses places. Autrement dit, c'est le dual tordu $\Hom(\,\wC_{F_\%},\mmu_{\ell^\%})$ du groupe $\,\wC_{F_\%}\simeq\Gal(F_\%^\lc/F_\%)$, où $\mmu_{\ell^\%}= \bigcup_n \mmu_{\ell^n}$ est le $\ell$-groupe des racines d'ordre $\ell$-primaire de l'unité.\smallskip

En résumé, il vient:

\begin{Th}\label{BP}
Soit $F$ un corps de nombres contenant les racines $2\ell$-ièmes de l'unité, $F_\%=\bigcup_n F_n$ sa $\Zl$-tour cyclotomique et $F_\%^\lc=\bigcup_n F_n^\lc$ la pro-$\ell$-extension abélienne maximale de $F_\%$ qui est partout complètement décomposée. Alors:\smallskip
\begin{itemize}
\item[(i)] Le groupe de Galois $\Gal(F_\%^\lc/F_\%)$ s'identifie à la limite projective $\,\wC_{F_\%}=\varprojlim \,\wCl_{F_n}$ des $\ell$-groupes de classes logarithmiques respectivement attachés aux étages finis $F_n$ de la tour.\smallskip

\item[(ii)] Le radical de Kummer $\Rad(F_\%^\lc/F_\%)$ s'identifie à la limite inductive $\,\wc_{F_\%}=\varinjlim \,\wc_{F_n}$ des pseudo-radicaux $\,\wc_{F_n}\simeq\,\T^\bp_{F_n}$, i.e. des modules de Bertrandias-Payan attachés aux $F_n$.\smallskip

\item[(iii)] Pour chaque $n\in\NN$, le $\ell$-groupe des classes logarithmiques de $F_n$ s'identifie au quotient de $\,\wC_{F_\%}$ fixé par $\,\Gamma_{\!n}=\Gal(F_\%/F_n)$; et le module $\,\T^\bp_{F_n}$ au sous-module de $\,\wc_{F_\%}$ fixé par $\Gamma_{\!n}$:\smallskip

\centerline{$\wCl_{F_n}\simeq {}^{\Gamma_{\!n}}\wC_{F_\%}$ \qquad \& \qquad $\T^\bp_{F_n} \simeq \,\wc_{F_n}= \,\wc_{F_\%}^{\,\Gamma_{\!n}}$.}
\end{itemize}
\end{Th}

\Remarque Le groupe $\,\wc_{F_n}$ {\em n'est pas}, en général, le radical de Kummer $\Rad(F_n^\lc/F_\%)$ attaché à la pro-$\ell$-extension abélienne maximale $F_n^\lc$ de $F_n$ qui est complètement décomposée partout sur $F_\%$.\smallskip

La descente kummérienne nécessite, en effet, une torsion {\em à la Tate}: si $\Tl=\varprojlim \mmu_{\ell^n}$ désigne le module de Tate construit sur les racines $\ell$-primaires de l'unité et $\lT = \varprojlim \mmu_{\ell^n}^*$ le module contagrédient, i.e. le dual de Pontryagin de $\mmu_{\ell^\%}$, on a:\smallskip

\centerline{$\Rad(F_n^\lc/F_\%) = (\lT \otimes_\Zl \wc_{F_\%})^{\Gamma_{\!n}}$ (et non $\,\wc_{F_\%}^{\,\Gamma_{\!n}}$).}

\newpage
\section{Rappels sur les éléments de Stickelberger normalisés}\label{2}

Soient $F$ un corps abélien de conducteur $f_F=f$ et $G_F=\Gal(F/\QQ)$ son groupe  de Galois.\smallskip

L'élément de Stickelberger normalisé attaché à $F$ est défini dans l'algèbre $\QQ[G_F]$ par la formule:\smallskip

\centerline{$\s_F^\ph=-\underset{0<a<f}{\sum^\times}\big(\frac{1}{2}-\frac{a}{f}\big)\;\big(\frac{F}{a}\big)^{\si{-1}}$,}\smallskip

\noindent où la somme porte sur les entiers $a$ étrangers à $f$ et $\big(\frac{F}{a}\big)$ désigne le symbole d'Artin. 
En particulier $\s_F^\ph$ n'est autre que la restriction à $F$ de l'élément $\s_f=\s_{\QQ[\zeta_f]}$ attaché au corps cyclotomique $\QQ[\zeta_f]$. 

Le symbole d'Artin $\sigma^\ph_a=\big(\frac{\QQ[\zeta_f]}{a}\big)$ est caractérisé par l'identité: $\zeta_f^{\sigma_a}=\zeta_f^a$. Et  $\sigma^\ph_{\si{-1}}=\big(\frac{\QQ[\zeta_f]}{-1}\big)$ est ainsi la conjugaison complexe, disons  $\bar\tau$.\smallskip
\begin{itemize}
\item Si donc $F$ est réel, on a $\big(\frac{F}{-1}\big)=1$, puis $\big(\frac{F}{f-a}\big)=\big(\frac{F}{-a}\big)=\big(\frac{F}{a}\big)$; et finalement $\s_F^\ph=0$.
\item Et, si $F$ est imaginaire, le calcul donne $\s_F^\ph=(1-\bar\tau)\,\s_F'$ avec $\s'_F=\underset{0<a<f/2}{\sum^\times}\big(\frac{1}{2}-\frac{a}{f}\big).\big(\frac{F}{a}\big)^{\si{-1}}$.
\end{itemize}
Il en résulte que l'élément de Stickelberger $\s_F^\ph$ est {\em imaginaire}, en ce sens qu'il se factorise par l'idempotent $e_-=\frac{1}{2}(1-\bar\tau)$ de l'algèbre $\QQ[G_F]$, l'idempotent complémentaire $e_+=\frac{1}{2}(1+\bar\tau)$ correspondant, lui, à la composante {\em réelle} de l'algèbre de Galois.\medskip

Considérons maintenant un sous-corps $K$ de $F$, autre que $\QQ$, de conducteur $f_K$ et de groupe de Galois $G_K$. Il est naturel de comparer $\s_K$ avec la restriction $N_{F/K}(\s_F^\ph)=\s^\ph_F|_K^\ph$ de $\s_F^\ph$ à $K$. Le résultat est le suivant (cf. e.g. \cite{Gra5} \S4.2):
\begin{equation}\label{Res}
N_{F/K}(\s_F^\ph)=\underset{p|f_F^\ph\;p\nmid f_K^\ph}{\prod}\,\big( 1-\big(\frac{K}{p}\big)^{\si{-1}}\big)\;\s_K^\ph
\end{equation}
où le produit porte sur les premiers $p$ qui divisent $f_F^\ph$ mais non $f_K^\ph$; autrement dit qui se ramifient dans $F/\QQ$ mais non dans $K/\QQ$. En particulier, on a toujours:\smallskip

\centerline{$N_{F/K}(\s_F^\ph)=\s_K^\ph$,}\smallskip

\noindent dès que $f_F^\ph$ et $f_K^\ph$ ont les mêmes facteurs premiers. En revanche, il vient:\smallskip

\centerline{$N_{F/K}(\s_F^\ph)=0$,}\smallskip

\noindent dès qu'il existe un premier $p$ ramifié dans $F/\QQ$ qui est complètement décomposé dans $K/\QQ$. C'est en particulier le cas lorsque $K$ est le sous-corps de décomposition d'un premier $p$ divisant $f$. Ainsi:

\begin{Prop}
L'élément de Stickelberger $\s_F^\ph$ annule les $\QQ[G_F]$-modules multiplicatifs engendrés par les idéaux premiers $\p_F^\ph$ de $F$ ramifiés dans $F/\QQ$. En d'autres termes on a: $\p_F^{\s_F}=1$ pour $\p_F | f$.
\end{Prop}

\Preuve On a, en effet: $\p^\ph_K=\p_F^{e_p(F/K)}$; puis: $\p_F^{e_p(F/K)\s_F }=\p_K^{N_{F/K}(\s_F)}=1$; d'où le résultat.
\medskip

Maintenant, pour obtenir des éléments entiers il est judicieux de tordre les éléments de Stickelberger par des multiplicateurs convenables; par exemple par les facteurs\smallskip

\centerline{$\delta_F^{\,c}= 1-c\big(\frac{F}{c}\big)^{\si{-1}}$,}\smallskip

\noindent où $c$ désigne un entier impair et étranger à $f$. Les éléments de Stickelberger tordus ainsi obtenus\smallskip

\centerline{$\s_F^{\,c}=\delta_F^{\,c}\,\s_F^\ph=\big( 1-c\big(\frac{F}{c}\big)^{\si{-1}}\big)\;\s_F^\ph$}\smallskip

\noindent sont alors dans $\ZZ[G_F]$. Un calcul immédiat donne, en effet:\smallskip

\centerline{$\s^{\,c}_F=\frac{1}{2}(c-1)\underset{0<p<f}{\sum^\times}\big(\frac{F}{a}\big)^{\si{-1}}\;-\;\frac{1}{f}\underset{0<p<f}{\sum^\times}\big(a\big(\frac{F}{a}\big)^{\si{-1}}-ac\big(\frac{F}{ac}\big)^{\si{-1}}\big)$;}\smallskip

\noindent avec à gauche un multiple de la norme $\nu_F^\ph=\sum_{\sigma\in G_F}\sigma$ et où la somme à droite est dans $f\ZZ[G_F]$.\medskip

Enfin, comme on a banalement $\delta_F^{\,c}|_K^\ph= 1-c\big(\frac{K}{c}\big)^{\si{-1}}=\delta_K^{\,c}$, ces éléments de Stickelberger tordus satisfont les mêmes identités de restriction que les éléments normalisés plus haut:
\begin{equation}\label{Res'}
N_{F/K}(\s_F^{\,c})=\underset{p|f_F^\ph\;p\nmid f_K^\ph}{\prod}\,\big( 1-\big(\frac{K}{p}\big)^{\si{-1}}\big)\;\s_K^{\,c}
\end{equation}
\noindent où $p$ parcourt les nombres premiers ramifiés dans $F/\QQ$ mais non dans $K/\QQ$.

\section{Dualité du miroir dans la $\Zl$-tour cyclotomique}

Supposons maintenant fixé un nombre premier impair $\ell$; partons d'un corps abélien  $F$  contenant le groupe $\mmu_\ell$ des racines $\ell$-ièmes de l'unité (ce qui implique en particulier que $\ell$ divise le conducteur $f$ de $F$); introduisons la $\Zl$-extension cyclotomique $F_\%=F\QQ_\%$ de $F$; et notons enfin $\Delta=\Gal(F_\%/\QQ_\%)\simeq\Gal(F/(F\cap\QQ_\%))$ le sous-groupe de $G_F=\Gal(F_\%/\QQ)$ qui fixe $\QQ_\%$.\smallskip

Écrivons $F_n$ le $n$-ième étage de la tour $F_\%$ (c'est un corps abélien imaginaire de conducteur $f_n=f\ell^n$); puis $G_n$ le groupe de Galois de $F_n/\QQ$ et $G_{F_\%}=\varprojlim G_n$ celui de $F_\%/\QQ$. 
Par projectivité, le groupe quotient $\Gal(\QQ_\%/\QQ)$, qui est isomorphe à $\Zl$, se relève dans $G_{F_\%}$, ce qui permet d'écrire $G_{F_\%}$ comme produit direct du sous-groupe fini $\Delta=\Gal(F_\%/\QQ_\%)$ et du relèvement $\Gamma=\gamma^{\Zl}\simeq \Zl$.\smallskip

L'algèbre de groupe complète $\Zl[[G_{F_\%}]]$ s'identifie alors à l'algèbre d'Iwasawa $\Zl[\Delta][[\gamma-1]]$ en l'indéterminée $\gamma-1$ à coefficients dans l'algèbre de groupe $\Zl[\Delta]$.

\begin{Prop}
Lorsque $\ell$ ne divise pas $|\Delta|$, l'algèbre $\Zl[\Delta]$ est un anneau semi-local, produit fini d'extensions non-ramifiées $\ZZ_\varphi$ de $\Zl$  indexées par les caractères $\ell$-adiques irréductibles $\varphi$ de $\Delta$:\smallskip

\centerline{$\ZZ_\varphi=\Zl[\Delta]e_\varphi$, avec $e_\varphi=\frac{1}{[\Delta |}\sum_{\sigma\in\Delta}\varphi(\sigma)\sigma^{\si{-1}}$,}\smallskip

Dans ce cas, l'algèbre complète $\Zl[[G_{F_\%}]]$ admet  la décomposition directe:\smallskip

\centerline{ $\Zl[[G_{F_\%}]]=\bigoplus_\varphi \Lambda_\varphi$,}\smallskip

\noindent  où  $\Lambda_\varphi =\ZZ_\varphi[[\gamma-1]]$ est l'algèbre d'Iwasawa attachée à l'extension non-ramifiée $\ZZ_\varphi$ de $\Zl$.\smallskip

Dans ce même contexte, le sous-corps  $F_0$ de $F^c$ fixé par $\Gamma$, linéairement disjoint de $\QQ^{\,c}$ sur $\QQ$ et de groupe de Galois $\Gal(F_0/\QQ)\simeq\Delta$, peut être pris dans $F$, qui est alors l'un des étages $F_m$ de la $\Zl$-extension cyclotomique  $F_0^{\,c}=F^c$ de $F^\ph_0$.
Et on a la relation entre conducteurs: $f=f_0\ell^m$.
\smallskip

Nous référons à cette situation particulière comme étant le cas semi-simple.
\end{Prop}

\Preuve La première partie de la Proposition est bien connue, la décomposition irréductible de l'algèbre $\Zl[\Delta]$ s'obtenant par relèvement de celle de l'algèbre semi-simple $\Fl[\Delta]$,  les facteurs simples (resp. locaux) de $\Fl[\Delta]$ (resp. de $\Zl[\Delta]$) correspondant aux caractères des représentations irréductibles de $\Delta$  sur $\Fl$ (resp. sur $\Zl$). L'exposant $e$ de $\Delta$ étant étranger à $\ell$, ces facteurs locaux sont des sous-extensions de $\Zl[\zeta_e]$ donc des extensions abéliennes non-ramifiées de $\Zl$.\par
Enfin, pour voir que $F$ provient par composition avec un étage fini $\QQ_m$ de la $\Zl$-tour $\QQ^{\,c}$ de $\QQ$, d'une extension abélienne $F_0$ de groupe $\Delta$, écrivons $G_{F_\%}=\Gamma\times\Delta$ et prenons un générateur $\gamma^{\ell^{\si{m}}}\delta$ du groupe procyclique $\Gal(F^c/F)$; notons $\ell^*\in\NN$ un inverse de $\ell$ modulo $|\Delta|$; et remplaçons $\gamma$ par $\gamma^\ph_0=\gamma\delta^{\ell^{\si{*\,m}}}$. Le sous-corps $F_0$ fixé par $\Gamma^\ph_0=\gamma_0^{\Zl}$ convient, $F$ étant, par hypothèse, fixé par $\Gamma_0^{\ell^{\si{m}}}$. \smallskip

Revenons maintenant au cas général:

\begin{Def}
Soit $\kappa : G_{F_\%} \to \ZZ_\ell^\times$ le caractère $\ell$-adique donné par l'action de $G_{F_\%}=\Gal(F_\%/\QQ)$  sur le groupe $\mmu^\ph_{\ell_\ph^{^{\si{\infty}}}}$ des racines de l'unité d'ordre $\ell$-primaire: $\zeta^\sigma=\zeta^{\kappa(\sigma)}\quad\forall\zeta\in\mmu^\ph_{\ell_\ph^{^{\si{\infty}}}}$.\smallskip

\noindent L'involution du miroir est définie sur l'algèbre de groupe $\Zl[[G_{F_\%}]]=\Zl[\Delta][[\gamma-1]]$ par l'identité:\smallskip

\centerline{$\Phi=\sum_{k\in\NN}a_k\,(\gamma-1)^k \mapsto \Phi^*=\sum_{k\in\NN}a^*_k\,(\kappa(\gamma)\gamma^{\si{-1}}-1)^k$,}\smallskip

\noindent où le reflet $a^*$ d'un élément $a=\sum_{\sigma\in\Delta}\alpha_\sigma\sigma\in\Zl[\Delta]$ est l'élément: $a^*=\sum_{\sigma\in\Delta}\alpha_\sigma\kappa(\sigma)\sigma^{\si{-1}}$.
\end{Def}

L'involution du miroir trouve son origine dans la dualité entre la description kummérienne des pro-$\ell$-extensions abéliennes d'un corps surcirculaire $F_\%$ et leur description galoisienne.\par
Le miroir envoie un élément $\rho$ de $G_{F_\%}$ sur l'élément $\rho^*=\kappa(\rho)\rho^{\si{-1}}$ (et inversement); d'où, par linéarité la formule annoncée sur $\Zl[G_{F_\%}]$ et, finalement, sur $\Zl[[G_{F_\%}]]=\Zl[\Delta][[\gamma-1]]$, la congruence $\kappa(\gamma) \equiv 1$ [mod $\ell$] assurant la convergence de la série.

\Nota La conjugaison complexe $\bar\tau$ opérant par passage à l'inverse sur les éléments de $\mmu^\ph_{\ell_\ph^{^{\si{\infty}}}}$, on a immédiatement $\kappa(\bar\tau)=-1$; donc $e_{\pm}^*=\frac{1}{2}(1\pm \bar\tau)^*=\frac{1}{2}(1\mp \bar\tau)=e_\mp$: le miroir échange composantes réelles et imaginaires.
Plus précisément, dans le cas semi-simple, il envoie l'idempotent $e_\varphi$ sur l'idempotent $e_{\chi\bar\varphi}$, où $\chi$ est la restriction de $\kappa$ à $\Delta$ et $\bar\varphi$ le contragrédient de $\varphi$: $\sigma\mapsto\varphi(\sigma^{\si{-1}})$.

\newpage
\section{Démonstration du Théorème Principal}

Supposons toujours $\ell$ premier impair et $F$ abélien contenant $\mmu_\ell^\ph$ de conducteur, disons, $f_F$; écrivons $F_\%$ la $\Zl$-extension cyclotomique de $F$ et $\Gamma=\gamma^\Zl$ le groupe procyclique $\Gal(F_\%/F)$.\par

Soient $\,\wCl_F$ le groupe des classes logarithmiques de $F$ et $\,\T^\bp_F$ le module de Bertrandias-Payan.

\begin{Prop}
Avec les notations précédentes et pour tout élément $c\in\Zl^\times$:
\begin{itemize}\smallskip

\item[(i)] Les éléments de Stickelberger tordus $\s_F^{\,c}=\s_{F_\%}^{\,c} \mod (\gamma-1)$ de $\Zl[G_F]$ annulent $\,\wCl_F$.\smallskip

\item[(ii)] Et leurs reflets respectifs $\s_F^{\,c*}=\s_{F_\%}^{\,c*} \mod (\gamma-1)$ annulent le groupe $\,\T^\bp_F$.
\end{itemize}
\end{Prop}

\Preuve Le théorème classique de Stickelberger (cf. e.g. \cite{Wa}, \S15.1), appliqué aux divers étages $F_n$ de la $\Zl$-tour cyclotomique $F_{\si{\infty}}/F$ affirme que, pour tout $c$ impair étranger aux conducteurs  $f^\ph_{F_n}$, les éléments de Stickelberger tordus $\s_{F_n}^{\,c}$ annulent respectivement les $\ell$-groupes de classes $\Cl_{F_n}$. Et comme les  $f^\ph_{F_n}$ ont les mêmes facteurs premiers que $f_F$, ces éléments vérifient les conditions de cohérence $N_{F_m/F_n}(\s_{F_m}^{\,c})=\s_{F_n}^{\,c}$ pour $m\ge n$, comme indiqué dans la section \ref{2}. Ainsi:\smallskip

$(i)$ Pour tout $c$ impair étranger à $f_F$ l'élément $\s_{F_\%}^{\,c}=\varprojlim \s_{F_n}^{\,c}$ de l'algèbre complète $\Zl[[G_{F_\%}]]$ annule le groupe $\C_{F_\%}\!=\varprojlim \Cl_{F_n}$, donc, en particulier, son quotient  $\,\C'_{F_\%}\!=\varprojlim \Cl'_{F_n}=\,\wC_{F_\%}$.\smallskip

Il suit de là que la classe $\s_F^{\,c}$ de $\s_{F_\%}^{\,c}$ modulo $(\gamma-1)$ annule le quotient ${}^\Gamma\wC_{F_\%}=\,\wCl_F$.\smallskip

$(ii)$ Par ailleurs, puisque $\s_{F_\%}^{\,c}$ annule $\,\wC_{F_\%}$, son reflet $\s_{F_\%}^{\,c*}$ annule $\,\wc_{F_\%}\simeq\Hom(\,\wC_{F_\%},\mmu_{\ell^\%})$. Et sa classe $\s_F^{\,c*}$ modulo modulo $(\gamma-1)$ annule donc le sous-module invariant $\wc_{F_\%}^{\,\Gamma}=\,\wc_F \simeq \,\T^\bp_F$.\smallskip

Faisant alors varier $c$, nous obtenons le résultat annoncé par un argument immédiat de densité.

\Remarque L'involution du miroir est définie dans l'algèbre complète $\Zl[[G_{F_\%}]]$. Ainsi $\s_F^{\,c*}$ désigne ici l'image dans $\Zl[G_F]=\Zl[[G_{F_\%}]]/(\gamma-1)$ du reflet $\s_{F_\%}^{\,c*}$ de l'élément $\s_{F_\%}^{\,c}$. Ce n'est pas {\em stricto sensu} le reflet de l'image $\s_{F}^{\,c}$ de $\s_{F_\%}^{\,c}$: pour obtenir une involution dans (un quotient de) $\Zl[G_F]$, il convient de raisonner modulo l'ordre $\ell^m$ de $\mu_F^\ph$, i.e. dans $(\ZZ/\ell^m\ZZ)[G_F]$ à l'instar de e.g. \cite{Gra5}.\medskip

Introduisons maintenant le symétrisé $\hat s^{\,c}_{F_\%}\!=s^{\,c}_{F_\%}\!+s^{\,c\,*}_{F_\%}\!$ de $s^{\,c}_{F_\%}$. L'élément $s^{\,c}_{F_\%}=\frac{1}{2}(1-\bar\tau)\hat s^{\,c}_{F_\%}$ est ainsi la partie imaginaire de $\hat s^{\,c}_{F_\%}$ et  $s^{\,c\,*}_{F_\%}=\frac{1}{2}(1+\bar\tau)\hat s^{\,c}_{F_\%}$ sa partie réelle. Et il vient:

\begin{Cor}
Pour chaque $n\in\NN$, la réduction $\hat s^{\,c}_{F_n}$ de $\hat s^{\,c}_{F_\%}$ {\rm mod} $\gamma^{\ell^n}-1$ annule le $\ell$-groupe $\,\wCl_{F_n}$ des classes logarithmique de $F_n$ ainsi que le module de Bertrandias-Payan $\,\T^\bp_{F_n}$.
\end{Cor}

\Preuve Notons $F_n^+$ le sous-corps réel de $F_n$ et $F^+_\%=F^+{}^z$ son unique $\Zl$-extension. De l'inclusion ${F_n^+}^\lc \subset {F_n^+}^\bp$, on conclut que $\,\wCl_{F_n^+}\simeq\Gal(F_n^+{}^\lc/F_\%^+)$ est un quotient de $\,\T^\bp_{F_n^+}\simeq\Gal(F_n^+{}^\bp/F_\%^+)$.\par

 La composante réelle $\,\wCl_{F_n}^+$ de $\,\wCl_{F_n}$ étant ainsi un quotient de celle de $\,\T^\bp_{F_n}$, l'élément $s^{\,c\,*}_{F_n}$ l'annule. Or, étant réel, il annule banalement la composante imaginaire $\,\wCl_{F_n}^-$; donc en fin de compte $\,\wCl_{F_n}$ tout entier. En résumé, $\hat s^{\,c}_{F_\%}$ annule $\,\wC_{F_\%}$ donc $\,\wc_{F_\%}$ par dualité; et $\hat s^{\,c}_{F_n}$ annule $\,\T^\bp_{F_n}$.\medskip

Venons-en enfin aux noyaux étales sauvages $\W_{2i}(F_n)$. Il est bien connu (cf. e.g. \cite{J44,Ng1}) qu'ils sont donnés, comme quotients des genres tordus, par les isomorphismes de modules galoisiens:\smallskip

\centerline{$\W_{2i}(F_n) \,\simeq\, {}^{\Gamma_{\!n}}(\,\Tl^{\otimes i}\otimes_\Zl \wC_{F_\%})=(\,\Tl^{\otimes i}\otimes_\Zl \wC_{F_\%})/(\,\Tl^{\otimes i}\otimes_\Zl \wC_{F_\%})^{(\gamma^{\ell^{\si{n}}}-1)}$,}\smallskip

\noindent où $\Tl^{\otimes i}$ est la $|i|$-ième puissance tensorielle de $\Tl$ pour $i > 0$, de $\lT$ pour $i < 0$ et $\Gamma_{\!n}=\gamma^{\ell^n\Zl}$.\smallskip

Introduisons les tordus {\em à la Tate} des symétrisés  $\hat s^{\,c}_{F_\%}\!$. Pour chaque élément $\Phi=\sum_{k\in\NN}a_k(\gamma-1)^k$ de l'algèbre de groupe $\Zl[[G_{F_\%}]]=\Zl[\Delta][[\gamma-1]]$ (avec $a=\sum_{\sigma\in\Delta}\alpha_\sigma\sigma\in\Zl[\Delta]$, $\Delta=\Gal(F_\%/\QQ_\%)$ et $\Gamma=\gamma^\Zl$ relevant $\Gal(\QQ_\%/\QQ)$), définissons le $i$-ième tordu {\em à la Tate} de $\Phi$ en posant:\smallskip

\centerline{$\Phi^{(i)}=\sum_{k\in\NN}a^{(i)}_k(\kappa(\gamma^i)\gamma-1)^k$, avec  $a^{(i)}=\sum_{\sigma\in\Delta}\alpha_\sigma\kappa(\sigma^i)\sigma$.}\smallskip

Avec ces conventions, il suit directement du corollaire précédent:

\begin{Cor}\label{Noyaux étales}
Pour chaque $n\in\NN$ et $i\in\ZZ$, la réduction $\hat s^{\,c\,(-i)}_{F_n}$  {\rm mod} $(\gamma^{\ell^n}-1)$ du $(-i)$-ième tordu à la Tate $\hat s^{\,c\,(-i)}_{F_\%}$ du symétrisé $\hat s^{\,c}_{F_\%}$ annule le $i$-ième $\ell$-noyau étale sauvage $\W_{2i}(F_n)$ attaché au $n$-ième étage $F_n$ de la tour cyclotomique $F_\%/F$.
\end{Cor}


\newpage

\def\refname{\normalsize{\sc  Références}}

\medskip\noindent
{\small
\begin{tabular}{l}
Institut de Mathématiques de Bordeaux \\
Université de {\sc Bordeaux} \& CNRS \\
351 cours de la libération\\
F-33405 {\sc Talence} Cedex\\
courriel : Jean-Francois.Jaulent@math.u-bordeaux.fr\\
{\footnotesize \url{https://www.math.u-bordeaux.fr/~jjaulent/}}
\end{tabular}
}


\addcontentsline{toc}{section}{Bibliographie}

{\footnotesize
\begin{thebibliography}{tt}

\bibitem{All} {\sc T. All},
{\em On $p$-adic annihilators of real ideal classes},
J. Number Theory {\bf 133} (2013), 2324--2338.


\bibitem{Ba} {\sc G. Banaszak},
{\em  Generalization of the Moore exact sequence and the wild kernel for higher $K$-groups},
Compos. Math. {\bf 86} (1993), 281--305.

\bibitem{BB} {\sc J. Barrett, D. Burns},
{\em Annihilating Selmer modules}, 
J. reine angew. Math. {\bf 675} (2013), 191--222.

\bibitem{BJ} {\sc K. Belabas, J.-F. Jaulent},
{\em The logarithmic class group package in PARI/GP},
Pub. Math. Besançon (2016).

\bibitem{BP} {\sc F. Bertrandias et  J.-J. Payan}, 
{\em  $\Gamma$-extensions et invariants cyclotomiques}, 
Ann. Sci. Éc. Norm. Sup. {\bf 4} (1972), 517--548.

\bibitem{CS} {\sc J. Coates,  W. Sinnott}, 
{\em An analogue of Stickelberger's theorem for the higher $K$-groups}, 
Invent. Math. {\bf 24} (1974), 149--161.

\bibitem{DK} {\sc S. Dasgupta, M. Kakde},
{\em On the Brumer-Stark conjecture},
 Prépublication (2020).

\bibitem{Gra1}{\sc G. Gras},
{\em Annulation du groupe des $\ell$-classes généralisées d'une extension abélienne réelle de degré premier},
Ann. Institut Fourier {\bf 29} (1979), 15--32.

\bibitem{Gra2}{\sc G. Gras},
{\em Class Field Theory: from theory to practice},
Springer Monographs in Mathematics (2005).

\bibitem{Gra3}{\sc G. Gras},
{\em Sur le module de Bertrandias-Payan dans une $p$-extension -- noyau de capitulation},
Pub. Math. Besançon  2016, 25--44.

\bibitem{Gra5}{\sc G. Gras},
{\em Annihilation of $tor_{\ZZ_p}(\mathcal G_{K,S}^\ab)$ for real abelian extensions $K/\QQ$},
Com. Adv. Math. Sci.{\bf 1} (2018), 5--34.
 
\bibitem{J23}  {\sc J.-F. Jaulent},
{\it  La Théorie de Kummer et le $K_2$ des corps de nombres},
J. Théor. Nombres Bordeaux {\bf 2} (1990), 377--411.
 
\bibitem{J28} {\sc J.-F. Jaulent}, 
{\em Classes logarithmiques des corps de nombres},
J. Théor. Nombres Bordeaux {\bf 6} (1994), 301--325.

\bibitem{J31} {\sc J.-F. Jaulent},
{\em Théorie $\ell$-adique globale du corps de classes},
 J. Théor.  Nombres Bordeaux {\bf 10} (1998),   355--397.

\bibitem{J55} {\sc J.-F. Jaulent},
{\em Sur les normes cyclotomiques et les conjectures de Leopoldt et de Gross-Kuz'min},
Annales Math. Québec {\bf 41} (2017), 119--140.

\bibitem{J56} {\sc J.-F. Jaulent},
{\em Sur la capitulation pour le module de Bertrandias-Payan},
Pub. Math. Besançon  2016, 45–58.

\bibitem{J66} {\sc J.-F. Jaulent},
{\em Annulateurs circulaires et conjecture de Grenberg},
Prépublication (2020), \url{https://arxiv.org/pdf/2003.12301 }

\bibitem{J44} {\sc J.-F. Jaulent, A. Michel},
{\em Approche logarithmique des noyaux étales des corps de nombres},
J. Number Th. {\bf 120} (2006), 72–91.

\bibitem{JN} {\sc  H. Johnston, A. Nickel}, 
{\em An unconditional proof of the abelian equivariant Iwasawa main conjecture and applications},
 Prépublication (2020).

\bibitem{Ng1} {\sc T. Nguyen Quang Do},
{\em  Analogues supérieurs du noyau sauvage},
J. Théor. Nombres Bordeaux {\bf 4} (1992), 263--271.

\bibitem{Ng2} {\sc T. Nguyen Quang Do},
{\em  Descente galoisienne et capitulation entre modules de Bertrandias-Payan},
Pub. Math. Besançon  2016, 59--79.

\bibitem{Ni} {\sc A. Nickel},
{\em  Annihilating wild kernels},
Doc. Math. {\bf 24} (2019), 2381--2422.

\bibitem{NV} {\sc T. Nguyen Quang Do, V. Nicolas},
{\em  Nombres de Weil, sommes de Gauss et annulateurs galoisiens}, 
Amer. J. Math. {\bf 133} (2011), 1533--1571.

\bibitem{Or1} {\sc B. Oriat},
{\em Annulation de groupes de classes réelles},
 Nagoya Math. J. {\bf 81} (1981), 45--56.

\bibitem{Sc} {\sc P. Schneider},
{\em  Über gewisse Galoiscohomologiegruppen},
Math. Z. {\bf 168} (1979), no. 2, 181--205.

\bibitem{Si} {\sc W. Sinnott}, 
{\em  On the Stickelberger ideal and the circular units of an abelian field}, 
Invent. Math. {\bf 62} (1980), 181--234.

\bibitem{Sn} {\sc P. Snaith},
{\em Relative $K_0$, annihilators, Fitting ideals and the Stickelberger phenomena},
Proc. London Math. Soc. {\bf 90} (1992), 545--590.

\bibitem{Sol} {\sc D. Solomon},
{\em On a construction of $p$-units in abelian fields}, 
Invent. Math. {\bf 109} (1992), 329--350.

\bibitem{Sou} {\sc C. Soulé},
{\em  $K$-théorie des anneaux d'entiers des corps de nombres et cohomologie étale},
Inv. Math. {\bf 55} (1979), 251--295.


\bibitem{Wa} {\sc L.C. Washington}, 
{\em Introduction to cyclotomic fields},
Graduate Texts in Math. {\bf 83} (1997).

\bibitem{We} {\sc C. Weibel},
{\em  The norm residue isomorphism theorem},
J. Topol. {\bf 2} (2009), 346--372.

\end{thebibliography}
}
 \end{document}